\documentclass[12pt,reqno,a4paper]{amsart}
\title[Global optimization in inverse problems]{ Global optimization in inverse problems:
A comparison of Kriging and radial basis functions }

\usepackage{amssymb}
\usepackage{setspace}

\def\epsilon{\varepsilon}
\def\phi{\varphi}
\def\theta{\vartheta}
\def\rho{\varrho}

\def\beq#1{
\begin{equation}\label{#1} {\tt\hspace*{-6em{[#1]~~~}}}
}
\def\eeq{\end{equation}}

\def\vsep#1{{\vrule height #1 depth 0pt width 0pt}}
\author{W. Jacquet, B. Truyen, P. de Groen, I. Lemahieu \& J. Cornelis}
\newtheorem{def1}{Definition}

\def\beq#1{\begin{equation}\label{#1}}
\def\eeq{\end{equation}}

\def\vsep#1{{\vrule height #1 depth 0pt width 0 pt}} 
\def\phi{\varphi}
\def\theta{\vartheta}

\begin{document}
\doublespacing

\maketitle
\begin{abstract}
We study global optimization (GOP) in the framework of non-linear inverse problems with a unique solution.
These problems are in general ill-posed.
Evaluation of the objective function is often expensive, as it implies the solution of a non-trivial forward problem. 
The ill-posedness of these problems calls for regularization while the high evaluation cost of the objective function can be addressed with response surface techniques.
The global optimization using Radial Basis Function (RBF) as presented by Gutmann  \cite{Gutmann2001} is a response surface global optimization technique with regularizing aspects.
Alternatively, several publications put forward global optimization using a probabilistic approach based upon Kriging as an efficient technique for
non-linear multi modal objective functions, thereby providing a credible stopping rule  \cite{Jones2001}. 
After comparing both concepts, we argue that in case of non-linear inverse problems an adaptation of the RBF algorithm seems to be the most promising approach.
\par\noindent \textbf{Keywords:} inverse problems,  Kriging, radial basis functions, global optimization (GOP)
\end{abstract}

\section{Introduction}
\setcounter{equation}{0}
The global optimization (GOP) technique using response surfaces developed in  \cite{Jones1998} and  \cite{Gutmann2001} addresses problems
where the only  information available is the ability to evaluate the objective function at high cost. 
The main idea is to replace the objective function by a response surface that interpolates the objective function over a finite set of evaluation points, for which its value is known.
The response surface is assumed to approximate the objective function and used as a surrogate for the objective function.
A naive idea is to search for the global minimum of the response surface, followed by the evaluation of the objective function at this minimum,
and the incorporation of this new data into the response surface,
after which the entire procedure is repeated using the updated response surface.
Unfortunately, the resulting sequence often does not contain a subsequence converging to a global minimum of the  objective function. 
A concentration of evaluation points can cause ``extrapolation'' type of errors in regions where the evaluation points are less dense.
Due to this phenomenon, the response surface can bridge an unexplored valley of the objective function and hide a possible global minimum there. 
A solution may consist of finding a good balance between minimum search and global exploration. 
Minimum search will result in a good approximation of the global minimum if the region of the global minimum is correctly represented by the response surface, and
if this surface does not introduce minima lower than the global minimum.
In order to guarantee a trustworthy representation, a more or less uniform coverage of the whole domain of interest should be generated.

\subsubsection*{Response Surface GOP, the Basic Algorithm:}
\begin{itemize}
\item \textbf{Choose} a number of \textbf{initial evaluation points}. \\
\textbf{Evaluate} the objective function at the evaluation points. \\
\textbf{Create} a \textbf{surrogate function} based upon the evaluation points and their evaluations,
called response surface, that should approximate the function to be minimized. \item
\begin{enumerate}
\item \textbf{Search for a new} \textbf{evaluation point}
      using the response surface and the position of the evaluation points,
      balancing scanning unexplored areas and believe in the response surface as a reliable surrogate.
\item \textbf{Evaluate} the objective function at this new evaluation  point, and update the surrogate function.
\end{enumerate}
\item \textbf{Repeat} $1$, $2$ until satisfaction or exhaustion.
\end{itemize}

This scheme can be used both with an RBF or Kriging approach.
The RBF algorithm and Kriging differ with regard to the interpolating function chosen and the way in which they balance exploration and minimum search.
In  what follows a concise description of standard RBF and the Kriging based GOP will be introduced followed by a comparison with non-linear inverse problems in mind.

\section{Response Surface GOP using Radial Basis Functions}

\setcounter{equation}{0}
The Radial Basis Function approach consists in constructing a linear space from which the interpolating functions are chosen dependent on the
arguments $x_1, \, \cdots \,,x_n$ - the position of the data points.
One way to introduce this dependency is to use a linear combination of functions with a radial symmetry about each of the $x_i$:
\begin{displaymath}
\sum_{i=1}^{n} \lambda_i \phi(\lVert x - x_{i} \rVert)
\end{displaymath}
where $\lVert \cdot \rVert$ is the Euclidean norm in $\mathbb{R}^d$ and $\phi$ is a function $[0,\infty[ \to
\mathbb{R}$.  \newline

Analytic results exist for the interpolation problem using this form, e.g.
the ``multiquadric'' $\phi(r)=  \sqrt{r^2+\gamma^2}$, $\gamma > 0$,
 and the ``Gaussian'' $\phi(r)=  e^{-\gamma r^2}$, $\gamma > 0$,  presented by Micchelli \cite{Micchelli1986}.
Despite these results it is useful to allow the addition of a polynomial of a low degree.
This extends the theory to include ``thin plate splines'', and to increases the number of
possible $\phi$'s to be used \cite{Powell1998}.  

\begin{def1}
Given $\{x_i \, \vert \, i=1 \, \cdots \, n\}$, $n$ pairwise distinct points in $\mathbb{R}^d$, an RBF interpolation function on $\mathbb{R}^d$ is a
function of the form :
\begin{displaymath}
s : \mathbb{R}^d \to \mathbb{R}: x \mapsto s(x)= \sum_{i=1}^{n}
\lambda_i \phi(\lVert x - x_{i} \rVert) + p(x)
\end{displaymath}
with $\lVert .  \rVert$ the Euclidean norm in $\mathbb{R}^d$,
 $\phi$ a function $[0,\infty[ \to \mathbb{R}$, and $p$ a polynomial of low degree with a limitation depending on the type of $\phi$.
\end{def1}
The functions $\phi$ considered by Jones \cite{Jones1998} are defined by one of the formulas in table 1.  Depending on the choice
for $\phi$ the RBF function, $s$ will be indicated with the adjective in the column ``RBF type''.
The degree of the polynomials will be restricted depending on the choice of $\phi$ related to the interpolating
capabilities of the corresponding RBF, as shown in the column ``Maximal degree'' of table 1, cf.  \cite{Powell1998}.
\begin{table}
$$
\begin{array}{|llc|c|c|}
\hline
         &                     & \,& \textrm{RBF type}          & \textrm{Maximal}        \\
         &                     & \,& \textrm{}                  & \textrm{degree}  \\
\hline \hline
\phi(r)= & r                   & \,& \textrm{linear}            & 0 \vsep{1.5em}\\
\phi(r)= & r^3                 & \,& \textrm{cubic}             & 1 \vsep{1.5em}\\
\phi(r)= & \left\{\begin{matrix} r^2log(r)  & {\rm if} & r > 0 \cr
                          0       & {\rm if} & r=0\end{matrix}\right.
                               & \,& \textrm{thin plate spline} & 1 \vsep{2.1em}\\
\phi(r)= & \sqrt{r^2+\gamma^2} & \,& \textrm{multiquadric}      & 0 \vsep{1.8em}\\
\phi(r)= & e^{-\gamma r^2}     & \,& \textrm{Gaussian}          & \textrm{no $p$ needed}\vsep{1.5em}\\
\hline
\end{array}
$$
\caption{Different choices for $\phi$, and the maximal degree $n$ of $p$.}
\end{table}

\subsection{Balance between exploration and minimum search}
Given the previously computed values of the objective function $f(x_1), \, \cdots \,,f(x_n)$, at the data points $x_1, \, \cdots \,,x_n$, the problem is to choose the next evaluation point. 
The idea is to start from an ``estimate'' of the value for the global minimum, then choose a ``reasonable'' candidate
for the location of the global minimum based upon a combination of this estimate, the response surface, and the position of the previously computed data points.
The main concern is limiting the introduction of possible artificial minima.
The choice is made such that within the class of interpolating functions a measure of ``bumpiness'' is minimal.

Let us have $n$ evaluation points of the objective function and its value there, $\{ (x_{i},f(x_{i}))  \, \vert  \, i = 1  \, \cdots \, n \}$.
Denote by $R_{n}$ the response surface that interpolates those points.
Assume that $f^{\ast}$ is some value below the calculated global minimum $R_{\min}:=\min \, R_n$ of the response surface $R_n$; $R_{\min}$ is an estimate of the global minimum of the objective function.
Obviously it is equal or below $f_{\min}:= \min_{i} f(x_i)$.
We now consider the response surface $\widetilde{R}_x$ that interpolates the objective function at the points $\{ (x_{i},f(x_{i}))  \, \vert  \, i = 1  \, \cdots \, n \} $, 
and takes the value $f^{\ast}$ at $x$.
This response surface $\widetilde{R}_x$ may be considered as a function of $x$.
The bumpiness of $\widetilde{R}_x$ obviously depends on the interpolating method and the location of the (fixed) evaluation points $\{x_i\}$.
It also depends on $x$ and the hypothetical value $f^{\ast}$.
If $f^{\ast}$ is near to $f(x_i)$ for some $x_i$, we may choose $x$ near $x_i$ without introducing steep gradients and without a large effect on the bumpiness of $\widetilde{R}_x$.
However, if $f^{\ast}$ is chosen far from any $f(x_i)$, we cannot choose $x$ near any $x_i$ if we do not want a large increase in bumpiness.
We can use this as a guidance in the selection of the next evaluation point.
If we constraint the search for the minimum by the condition that the increase of bumpiness of $\widetilde{R}_x$ is as small as possible,
then a choice of $f^{\ast}$ slightly below the calculated estimate of the minimum of $R_n$ is likely to result in an $x$ near the location of this minimum.
On the other hand, if the target value $f^{\ast}$ is far below, this may result in a point far away, in a yet unexplored part of the domain of interest.

\subsubsection*{Response Surface GOP using Radial Basis Functions, the Basic Algorithm:}
\begin{itemize}
\item \textbf{Choose} a number of \textbf{initial evaluation points}. \\
\textbf{Evaluate} the objective function at the evaluation points. \\
\textbf{Create} a \textbf{surrogate function}: determine the RBF $R_n$ interpolating the objective function at the evaluation points.
\item
\begin{enumerate}
\item \textbf{Search for a new} \textbf{evaluation point}:
\begin{itemize}
\item determine the global minimum of the surrogate function $R_n$ and decide upon a target value $T$ lower than its minimum,
\item      find a point $x^{\star}$ that minimizes the chosen measure of bumpiness of the RBF interpolating $\{ (x_{i},f(x_{i}))  \, \vert  \, i = 1  \, \cdots \, n \} \cup \{ (x,T)  \}$.
\end{itemize}
\item \textbf{Evaluate} the objective function at this new evaluation  point, and update the surrogate function.
\end{enumerate}
\item \textbf{Repeat} $1$, $2$ until satisfaction or exhaustion.
\end{itemize}

\section{Response Surface GOP using Kriging}
\setcounter{equation}{0}
A wide research field emerged around the usage of probabilistic models, such as Kriging, in the field of global optimization.
A concise description of the use of Kriging in the context of global optimization can be found in \cite{Jones2001},
whereas a more elaborate treatment is given in \cite{Jones1998}.
In these papers, the  validity is claimed using the statistical terminology and frequency interpretations.
The statistical basis for spatial data analysis can be found in  \cite{Anderson1984} and \cite{Cressie1993}.

Kriging starts from the fundamental hypothesis that the objective function to be optimized is the realization of a random field $Y$.
For convenience , it is assumed that this random field is square integrable.
As a consequence, the  modelling of the random field can be split into the modelling of 
a random field with $0$ expectation at each point of its domain, and the modelling of the expectation. \newline

The only information about the objective function resides in the calculated values at the observation points. 
For each random variable $Y_{x_i}$ located at an observation point $x_i$ only one observation is available: the calculated value $f(x_i)$ of the objective function at this point.
Meaningful statistics, however, are difficult to extract starting from just one value.
To overcome this problem additional model hypotheses are needed, linking the random variables of the different observation points. 
One approach is to adopt a form of ``stationarity'' hypothesis.
``Stationarity'' stems from the continuous-time series environment, meaning that the probabilistic structure in some sense does not depend on a shift in time, and
the covariance between the variables at different moments in time merely depends on the time lag.  
In the analysis of spatial data this comes down to the assumption that the probabilistic structure does not depend on the location
 and the covariance between two variables at different places in one direction merely depends on the distance between them. \newline

It is assumed that the mean  and the variance of the variables of the random field are constant for all locations $u$ in the region of interest $D$ of the GOP,
$E[Y(u)]=\mu$, $\textrm{var} ( Y(u)) = \sigma^2$.
The next step is to introduce a parametrized model for the covariance structure of the random field:\newline

\begin{equation}
Corr(Y(u),Y(v)) = \exp( - \sum_{j=1}^{d} \theta_j \lvert u_j - v_j \rvert^{p_j})
\end{equation}
with $u$ and $v$ points in the acceptance region $D$.  
The last step is to introduce a probability density hypothesis:\newline
For all $u_1,  \, \cdots \, , u_m$ in $D$ the multivariate random variable $(Y(u_1), \, \cdots \,,Y(u_m))$ follows a multivariate normal distribution.  
The parameters $\theta_1 ,  \, \cdots \, , \theta_d \ge 0$, $p_1 ,  \, \cdots \, , p_d \ge 0$, $\mu$ and $\sigma^2$ are to be estimated using the partial observation.
This estimation can be done using a ``two step'' loglikelihood estimation technique described in \cite{Jones2001} and \cite{Jones1998}.  
The estimation in itself is a non-trivial global optimization problem.
Starting from the parameters a prediction can be made for the unobserved part of the realization of the random field.
This results in a function $s$ that interpolates the partial realization of the following form:
\begin{displaymath}
s(x)= \sum_{i=1}^{n} \, \widehat{\lambda}_i \,  \exp \left(
{\textstyle - \sum_{j=1}^{d} \widehat{\theta}_j \lvert x - x_{i}
\rvert^{\widehat{p}_j} } \right) + \widehat{p}
\end{displaymath}
The Kriging interpolating function $s$ is a linear combination of radial basis functions and a constant.

\subsection{Balance between exploration and minimum search}

Together with an estimate of the objective function, which is a realization of the random field, several auxiliary functions can be estimated at each point in the acceptance region such as
the lower and upper confidence interval functions, the probability of improvement, and the expected improvement.
The probability of improvement and expected improvement are suited for  finding  a new observation point.
The probability of improvement, given a target value $T$,  $P[Y(x) \le T]$, amounts to the probability  that the value at a certain point will be lower than the chosen target $T$.
The expected improvement at a point $x$, $E[Y(x)-f_{\min}]$, is the expectation  of the difference between the random field $Y(x)$,  and the constant value $f_{\min}$. 
This is being used in a GOP algorithm, choosing target values $T$ alternating between values significantly smaller than $f_{\min}$ and values closely approaching $f_{\min}$.
Maximizing the probability of improvement for a well chosen cycle of $T$, balances scanning unexplored areas and believe in the response surface
as an approximation of the objective function.
The result of maximizing the expected improvement does not require a target value to be set.

\subsubsection*{Response Surface GOP using Kriging, the Basic Algorithm:}
\begin{itemize}
\item \textbf{Choose} a number of \textbf{initial evaluation points}.\\
\textbf{Evaluate} the objective function at the evaluation points.\\
\textbf{Create} a \textbf{surrogate function} based upon the evaluation points and their evaluations:
\begin{itemize}
\item estimate the scaling and smoothness parameters minimizing the maximum likelihood function
\item calculate the coefficients of the Kriging interpolating function
\end{itemize}
\item
\begin{enumerate}
\item \textbf{Search for a new} \textbf{evaluation point}:
\begin{itemize}
\item      determine $f_{\min}$ the minimum of the currently known values  and decide upon a target value $T$ lower than $f_{\min}$,
           find a point $x^{\star}$ that maximizes the probability of improvement with respect to $T$,
\item      or find a point $x^{\star}$ that maximizes the expected improvement $E[Y(x)-f_{\min}]$.
\end{itemize}
\item \textbf{Evaluate} the objective function at this new evaluation  point, and update the surrogate function.
\end{enumerate}
\item \textbf{Repeat} $1$ till $2$ until satisfaction or exhaustion.
\end{itemize}

\section{Discussion and Comparison of RBF and Kriging GOP in Inverse Problems}
\setcounter{equation}{0}
Kriging replaces the linear combination of radially symmetric functions about the observation points by a linear combination of functions with ``directional'' symmetry about the observation points.
When the parameters of the Kriging model are restricted by stating that all parameters $\theta_i$ as well as $p_i$ are equal, the interpolating estimate becomes a RBF. 
Even with these restrictions, the Kriging algorithm does not become an RBF GOP.
Indeed the algorithms still differ in the criterion used in the selection of a new observation point.
The search for the new observation point in the RBF GOP stems from an energy reasoning, and the direct need to have an automatic equilibrium between exploration and believe.
In the Kriging GOP on the other hand this equilibrium is a consequence of  criteria developed using probabilistic reasoning.
The basic assumption underlying the development of this probabilistic theory is that the observations are part of one realization of a random field.
In a large class of inverse problems, the objective function is deterministic.
Although it is not hard to introduce a random field with a given objective function as possible outcome,
it is however, difficult to interpret probabilistic statements in a meaningful way if the random field is not theoretically linked to the problem at hand.
Even if  the objective function can be interpreted as the outcome of a random field, a second more fundamental problem related to the probabilistic theory emerges.
In the iterative process of finding a global optimum, the observation points are chosen so as to obtain a sequence of values for the objective
function that converges to its global minimum.
In this way, the new observation point is chosen based upon all previous observations.
Nevertheless, in the development of the parameter estimators, the Kriging approach treats the observation points and their corresponding function values 
 as if the position of one observation point in the sequence does not depend on the values of the objective function in the previous observation points.
If one tries to find a global minimum, it can be expected that most of the evaluations are below the ``average'' of the complete but unknown realization (the objective function).
In this way the estimate  of the mean will be an underestimation.
When the sequence of function values converges to the global minimum, the estimate of the variance will be an underestimation too.

The Kriging procedure, implies that the basis functions are given an orientational preference,
 which appears only to be a good idea if at every local minimum the eigenvalues and eigenvectors of the
 Hessian matrix of the objective function are similar.
When scouting about a potential local minimum the behavior of the objective function in the neighborhood of this minimum will tend to dominate the estimates.
If there is no clear evidence that this behavior can be generalized to all points in the acceptance region, certain regions will tend to be less favorably approximated
than in an approach with no orientational adaptation. 
In geological applications \cite{Cressie1993}, this re-scaling option agrees with the fact that the ``depth'' coordinate defined by the gravitational force is
to be treated distinctly from the two other coordinates on the surface, taken perpendicular to the depth coordinate.
However, this does not necessarily apply to other applications.

What remains is that in the Kriging GOP approach each of the criteria to chose the new evaluation point, seeks to find a balance between exploration and believe in the response surface as a surrogate.
Under certain conditions minimization of the expected improvement creates a subset, dense in the acceptance region, and thus will come close to the global minimum \cite{Locatelli1997}. 
These properties remain also true outside the Kriging setting, independent of the statistical validity arguments. 
There is no need to confine the criteria used in the recruitment of the new evaluation points, exclusively to the Kriging GOP setting.
Instead, it can be useful to examine their performance also in a RBF GOP setting. 

The RBF GOP algorithm replaces the original GOP problem by a sequence of GOP problems with a relatively simple update from one iteration to the next, cf. \cite{Bjorkman2000}.
The Kriging GOP algorithm instead replaces the original GOP problem by a sequence of two GOP problems, the first being the non-trivial estimation of the model parameters.

In inverse problems there is a clear need for regularization.
When working with the Gaussian RBF's, this can be done by controlling the impact of the radial basis function via 
the parameter $\gamma$.
Following the Kriging algorithm the adaptation of the parameters to the observations might cause the
local behavior of the response surface to be modelled on a too detailed level,
cancelling out the regularization effect necessary to compute a meaningful solution to ill posed inverse problems.

\section{Conclusions}
\setcounter{equation}{0}

The application of a  probabilistic model to find the solution of an inverse problem
is not evident, both from theoretical and practical point of view.
Not only is there a clear cost in the adaptation of the basis functions to the observations,
also the adaptation itself is not necessarily  wanted,
 because it may cancel the underlying regularizing effects, and might model the ``noise''.
Instead it could be interesting to use adaptive RBF's, with automatic updating based upon the data, incorporating the fact that the sequence of observation points is not random.
When doing so, safeguards have to be build in to retain the necessary regularization properties of the method.
The criterion used in the selection of the new iteration point in the Kriging algorithm, may be adapted to the RBF GOP algorithm.

\bibliographystyle{plain}
\bibliography{Kriging_GOP}

\newpage

Wolfgang Jacquet
Department of Electronics and Information Processing - ETRO,
Vrije Universiteit Brussel,
Pleinlaan 2,
B-1050 Brussels,
Belgium

Bart Truyen
Department of Electronics and Information Processing - ETRO,
Vrije Universiteit Brussel,
Pleinlaan 2,
B-1050 Brussels,
Belgium

Pieter de Groen
Department of Mathematics,
Vrije Universiteit Brussel,
Pleinlaan 2,
B-1050 Brussels,
Belgium

Ignace Lemahieu
ELIS/MEDISIP,
UGent,
St.-Pietersnieuwstraat 41,
B-9000 Gent,
Belgium

Jan Cornelis
Department of Electronics and Information Processing - ETRO,
Vrije Universiteit Brussel,
Pleinlaan 2,
B-1050 Brussels,
Belgium

This work was supported in part by the FWO (Fund for Scientific Research - Flanders (Belgium)) project G.0072.01,
 "Stabilized deconvolution methods for inverse problems with application to linear (Magnetic Resonance imaging) and nonlinear (Ground Penetrating Radar imaging) image reconstruction,´´
 and the Concerted Research Action (GOA) of the Research Council of the Vrije Universiteit Brussel "Numerical issues in tomographic shallow subsurface imaging with application to land mine detection´´

\end{document}